\documentclass[reqno,A4paper]{amsart}
\usepackage{mathrsfs,mathtools,latexsym,eucal}
\usepackage{amscd,amsfonts,amssymb,amsmath}
\usepackage{graphicx,graphics,color}
\usepackage{geometry,setspace}
\usepackage[utf8]{inputenc}
\usepackage[active]{srcltx}
\usepackage[english]{babel}
\usepackage[pagewise,displaymath,mathlines]{lineno}
\usepackage{enumerate}
\usepackage{dsfont}
\newtheorem{theorem}{Theorem}[section]
\newtheorem{lemma}[theorem]{Lemma}
\newtheorem{proposition}[theorem]{Proposition}
\newtheorem{corollary}[theorem]{Corollary}

\newtheorem{definition}[theorem]{Definition}
\newtheorem{example}[theorem]{Example}

\numberwithin{equation}{section}
\newcommand{\N}{\mathbb{N}}

\newcommand{\R}{\mathbb{R}}

\newcommand{\disp}{\displaystyle}
\newcommand{\vareps}{\varepsilon}

\def\nequiv{\hspace{.17cm} / \kern-.99em \equiv }
\def\IC{\mathbb{C}}
\def\IN{\mathbb{N}}

\def\IR{\mathbb{R}}

\def\lq{\leqslant}
\def\gq{\geqslant}

\begin{document}
\title{Measures on effect algebras}
\author{G. Barbieri, F. J. Garc\'ia-Pacheco, S. Moreno-Pulido}
\newcommand{\acr}{\newline\indent}
\address{\llap{*\,}Department of Mathematics\acr
                   University of Salerno\acr
                   Via Giovanni Paolo II\acr
                   84084 Fisciano\acr
                   ITALY}
\email{gibarbieri@unisa.it}

\address{\llap{**\,}Department of Mathematics\acr
                    College of Engineering\acr
                    University of Cadiz\acr
                    Puerto Real 11519\acr
                    SPAIN}
\email{garcia.pacheco@uca.es, soledad.moreno@uca.es}




\begin{abstract}
We study measures defined on effect algebras. We characterize real-valued measures on effect algebras and find a class of effect algebras, that include the natural effect algebras of sets, on which $\sigma$-additive measures with values in a finite dimensional Banach space are always bounded. We also prove that in effect algebras the Nikodym and the Grothendieck properties together imply the Vitali-Hahn-Saks property, and find an example of an effect algebra verifying the Vitali-Hahn-Saks property but failing to have the Nikodym property. Finally, we define the concept of variation for vector measures on effect algebras proving that in effect algebras verifying the Riesz Decomposition Property, the variation of a finitely additive vector measure is a finitely additive positive measure.
\end{abstract}

\maketitle
\section{Introduction}

The concept of effect algebra was originally introduced in \cite{FouBen94}. An effect algebra is a four-tupla $(L,\oplus,0,1)$ consisting of a set $L$ with two special elements $0,1\in L$, called the \emph{zero} and the \emph{unit}, and with a partially defined binary operation $\oplus$ satisfying the following conditions for all $p,q,r\in L$:
\begin{enumerate}
\item[(E1)] \emph{Commutativity:} If $p\oplus q$ is defined, then $q\oplus p$ is defined and $p\oplus q= q\oplus p$.
\item[(E2)]  \emph{Associativity:} If $q\oplus r$ is defined and $p\oplus (q\oplus r)$ is defined , then $p\oplus q$ and $(p\oplus q)\oplus r$ are defined and $(p\oplus q)\oplus r=p\oplus(q\oplus r)$.
\item[(E3)]  \emph{Orthosupplementation:} For every $p\in L$, there exists a unique $q\in L$, called the orthosupplement of $p$, such that $p\oplus q=1$.
\item[(E4)]  \emph{Zero-One Law:} If $1\oplus p$ is defined , then $p=0$.
\end{enumerate}

A partial order can be defined in an effect algebra $L$: $p\leq q \Leftrightarrow \exists r\text{ with } p\oplus r = q.$ Every effect algebra obeys the cancellation law and in every effect algebra a dual operation $\ominus$ to $\oplus$ can be defined as follows: $a\ominus c $ exists and equals $b$ if and only if $b\oplus c$ is defined and $b\oplus c=a$. With this notation the orthocomplement of an element $a$ is $1\ominus a$. Whenever we write $p\perp q$ we will be assuming that $p\oplus q$ is defined and we will say that $p$ and $q$ are orthogonal. The orthocomplement of an element $a$ is also denoted by $a^\perp$.\\

If $a_1,\dots, a_n\in L$, we inductively define $a_1\oplus\dots\oplus a_n=(a_1\oplus\dots\oplus a_{n-1})\oplus a_n$ provided that the right-hand side exists. The definition is independent on permutations of the elements. We say that a finite subset $\{a_1,\dots ,a_n\}$ of $L$ is orthogonal if $a_1\oplus\dots\oplus a_n$ exists. A subset $G\subset L$ is said to be orthogonal if so is every finite subset of $G$. Furthermore, if $\sup\{\oplus F:F\subset G,F\mbox{ finite}\}$ exists, then we define $\oplus G:=\sup\{\oplus F:F\subseteq G, F\mbox{ finite}\}$. An effect algebra in which every orthogonal sequence $(a_n)_{n\in\N}$ verifies that $\bigoplus_{n\in\N}a_n$ exists is called a $\sigma$-effect algebra.\\

Particular examples of effect algebras are the so called effect algebras of sets, that is, subsets $L$ of the power set $\mathcal{P}(X)$ of a given set $X$ such that $\{\varnothing,X\}\subseteq L$ and $(L,\cup,\varnothing,X)$ has effect algebra structure under the partial operation $A\cup B$ defined only for those $A,B\in L$ such that $A\cap B=\varnothing$. A natural effect algebra of sets is an effect algebra of subsets of $\N$ containing the finite subsets of $\N$.\\

An effect algebra $L$ is said to enjoy the Riesz Decomposition Property (RDP) (see \cite{DvuPulma00}) if $c\leq a \oplus b$ implies that $c=c_1\oplus c_2$ with $c_1\leq a$ and $c_2\leq b$. Notice that this condition easily extends to finite orthosums (see Lemma \ref{rdpl}). If it also holds for infinite orthosums, then we say that the effect algebra verifies the strong Riesz Decomposition Property (sRDP). In other words, (sRDP) means that if $c\leq \bigoplus_{n\in\N}a_n$, then $c=\bigoplus_{n\in\N}c_n$ and $c_n\leq a_n$ for all $n\in\N$. The standard scale effect algebra $L:=\IR^+[0,1]=\{x\in\IR:0\lq x\lq 1\}$ where $x\oplus y=x+y \mbox{ if and only if } x+y\lq 1,$ verifies the RDP.\\

We will finish the introduction with the definition of measure on an effect algebra. As far as we know, while modular measures are extensively studied among others by Avallone, Barbieri, Basile,  Vitolo, Weber (e.g., see \cite{Av, Ab, B, B1, W}), measures on effect algebras are not studied so much. For some interesting properties about measures we refer the reader to \cite{AiMoRa13, AiTa, HwangHongRonglu, JunZhi}.\\

Here, measures on effect algebras are valued in (possibly normed) groups (all groups considered throughout this manuscript will be commutative and additive). So, let us recall first the definition of norm for a group. A semi-norm on a group $G$ is a function $\|\cdot\|:G\to[0,\infty)$ verifying that:
\begin{enumerate}
\item Identity preserving: $\|0\|=0$.
\item Symmetricity: $\|-g\|=\|g\|$ for all $g\in G$.
\item Triangular inequality or sub-additivity: $\|g+f\|\leq \|g\|+\|f\|$ for all $g,f\in M$.
\end{enumerate}
If, in addition, it is verified that $\|g\|=0$ implies $g=0$, then the semi-norm is called a norm.


\begin{definition}\label{mea}
Let $L$ be an effect algebra, $G$ a topological (or normed if necessary) group and $\mu:L\rightarrow G$ a map.
\begin{itemize}
\item $\mu$ is a \emph{measure} if it is additive, that is, for every orthogonal $a,b\in L$, $\mu(a\oplus b)=\mu(a)+\mu(b)$.
\item $\mu$ is a \emph{countably additive or $\sigma$-additive measure} if for every orthogonal sequence ${(a_n)}_{n\in\N}\subseteq L$ such that $\bigoplus_{n\in\IN} a_{n}$ exists, $\mu\left(\bigoplus_{n\in\IN} a_{n}\right)=\sum_{n\in\IN}\mu(a_n)$.
\item $\mu$ is a \emph{strongly additive measure} if for every orthogonal sequence ${(a_n)}_{n\in\N}\subseteq L$, $\sum_{n\in\IN}\mu(a_n)$ converges.
\item $\mu$ is an \emph{absolutely additive measure} if for every orthogonal sequence ${(a_n)}_{n\in\N}\subseteq L$, \break$\sum_{n\in\IN}\|\mu(a_n)\|<+\infty$.
\item $\mu$ is a \emph{bounded measure} if $\sup\{\|\mu(a)\|:a\in L\}<+\infty$.
\item $\mu$ is a \emph{exhaustive or strongly bounded measure} if for every orthogonal sequence ${(a_n)}_{n\in\N}\subseteq L$, $\lim \mu(a_n)=0$.
\end{itemize}\end{definition}

\section{Characterization of measures on effect algebras}

It is well known that the last four items of Definition \ref{mea} are in fact equivalent for real measures defined on Boolean algebras (see \cite{DiUhl77}). We will show next that for real measures on effect algebras all of them are equivalent except for the boundedness.

\begin{theorem}\label{equiv}
Let $L$ be an effect algebra, $G$ be a complete Hausdorff topological group and $\mu:L\rightarrow  G$ a measure.
\begin{enumerate}
    \item Then $\mu$ is strongly additive if and only if it is exhaustive.

    \item If $G=\R$, then $\mu$ is absolutely additive if and only if $\mu$ is strongly additive. Moreover, if $\mu$ is bounded, then $\mu$ is exhaustive.
\end{enumerate}
\end{theorem}
\begin{proof}
\mbox{}
\begin{enumerate}

\item We reduce to the classical case. Let $(a_n)_{n\in\N}$ be an orthogonal sequence in $L$, and observe that $\bar\mu:\phi_0(\N)\rightarrow G$ defined by $\bar\mu(F):=\mu(\oplus_{n\in F}\;a_n)$ is a measure on the ring of sets $\phi_0(\N):=\{A\subset \N: A\text{ is finite}\}$ which is exhaustive if and only if it is strongly additive (see \cite[Corollary 18, p. 8]{DiUhl77}).

\item On the first hand, the equivalence holds true since unconditionally convergence is equivalent to absolute convergence in real series. Finally, we will show that if $\mu$ is bounded, then it is exhaustive. Suppose that $\mu$ is bounded, but not exhaustive and let $(a_n)_{n\in\N}\subseteq L$ be an orthogonal sequence such that $a_n\nrightarrow0$. Then  $\sum_{n=1}^{\infty}|\mu(a_n)|=+\infty$. Consider the sets $M^+ =\{n\in\IN:\mu(a_n)>0\}$, $M^-=\{n\in\IN:\mu(a_n)<0\}$. Then $\sum_{n\in M^+}\mu(a_n)=+\infty$ or $\sum_{n\in M^-}\mu(a_n)=-\infty$. Suppose (without lack of generality) that $\sum_{n\in M^+}\mu(a_n)=+\infty$, this is, $\sum_{k=1}^\infty\mu(a_{n_k})=+\infty$, with $\mu(a_{n_k})>0$, for all $k\in\IN$. Then
$$\sum_{k=1}^\infty\mu(a_{n_k})=\lim_{m\to+\infty}\sum_{k=1}^m\mu(a_{n_k})=\lim_{m\to+\infty}\mu\biggl(\bigoplus_{k=1}^m a_{n_k}\biggr)=+\infty,$$
which contradicts the fact that $\mu$ is bounded.

\end{enumerate}

\end{proof}


We will next present an example (Example \ref{ejemploEfectoNoAcotadaSumable}) of an unbounded absolutely additive measure. However, we will first demonstrate a lemma. Recall that $\phi(\IN):=\{A\in\mathcal{P}(\IN):A\mbox{ is finite or cofinite}\}$.

\begin{lemma}\label{ExistenciaBn}
There exists a sequence ${(B_k)}_{k\in\IN}$ of subsets of $\IN$ such that
\begin{enumerate}[(1)]
\item $B_i\cap B_j\notin\phi(\IN)$.
\item $B_i\cap B_j^c\notin\phi(\IN)$, if $i\neq j$.
\item $B_i^c\cap B_j^c\notin\phi(\IN)$.
\end{enumerate}\end{lemma}
\begin{proof}
Denote by $P=\{p_n: n\in\IN\}$ the sequence of prime numbers of $\IN$ and consider $A_n:=\{p_n^m:m\in\IN\}$ for every $n\in\IN$. Note that $A_n\cap A_m=\varnothing$ if $n\neq m$ and $A_n\notin\phi(\IN)$ for every $n\in\IN$. Define now $$B_k:=\bigcup_{n=1}^k A_{2n-1}\cup\bigcup_{n=k}^\infty A_{2n}$$ for $k\in\IN$. Notice that $B_k\notin\phi(\IN)$ for every $k\in\N$. These sets verify the desired properties:
\begin{enumerate}
\item $B_i\cap B_j\notin\phi(\IN)$. Indeed, $(B_i\cap B_j)^c=B_i^c\cup B_j^c$ is a union of infinite sets so it is clearly infinite; and $A_1\subset B_i\cap B_j$ so $B_i\cap B_j$ is also infinite.
\item $B_i\cap B_j^c\notin\phi(\IN)$. Indeed, $(B_i\cap B_j^c)^c=B_i^c\cup B_j$ is again infinite for being a union of infinite sets; and $A_{2i}\subseteq B_i\cap B_j^c$ if $i< j$ and $A_{2i-1}\subseteq B_i\cap B_j^c$ if $i>j$, so $B_i\cap B_j^c$ is also infinite for $i\neq j$.
\item $B_i^c\cap B_j^c\notin\phi(\IN)$. Indeed, $(B_i^c\cap B_j^c)^c=B_i\cup B_j$ is again infinite for being a union of infinite sets; and $A_{2(i+j)+1}\subseteq B_i^c\cap B_j^c$ so $B_i^c\cap B_j^c$ is also infinite.
\end{enumerate}
\end{proof}

We will show now an example in which we illustrate that in the context of effect algebras a real measure can be absolutely additive but not bounded.

\begin{example}\label{ejemploEfectoNoAcotadaSumable}
According to Lemma \ref{ExistenciaBn}, consider a sequence ${(B_n)}_{n\in\IN}$ of subsets of $\IN$ such that:
\begin{itemize}
\item $B_i\cap B_j\notin\phi(\IN)$.
\item $B_i\cap B_j^c\notin\phi(\IN)$, if $i\neq j$.
\item $B_i^c\cap B_j^c\notin\phi(\IN)$.
\end{itemize}
Let $L:=\{\varnothing\}\cup\{\IN\}\cup\{ B_n:n\in\IN\}\cup\{B_n^c:n\in\IN\}$ and $\mu:L\rightarrow\IR$ such that $\mu(\varnothing) =\mu(\IN)=0$, $\mu(B_n)=n$ and $\mu(B_n^c) =-n$. Observe that $(L,\cup,\varnothing,\IN)$ is an effect algebra, where $\cup$ is defined by $A\cup B\in L\Leftrightarrow A,B\in L$ and $A\cap B=\varnothing$. Note that by the definition of the $B_n$'s, the orthogonal sequences ${(E_n)}_{n\in\IN}$ in $L$ are trivial, that is, there exists $n_0\in\IN$ such that $E_n=\varnothing$ for every $n\gq n_0$ and thus $\disp\sum_{n=1}^\infty|\mu(E_n)|<+\infty$. However, $\mu$ is not bounded, because $\mu(B_n)=n$ for every $n\in\IN$.\end{example}

The previous example showed a real measure on a effect algebra which is absolutely additive but not bounded. We will show now that, on certain effect algebras which include the natural effect algebras of sets, $\sigma$-additive measures with values in a finite dimensional Banach space are necessarily bounded.

\begin{definition}
Let $L$ be a poset. A generator of $L$ is a subset $B$ of $L$ with the property that for every $a\in L$ there exists an increasing $(b_n)_{n\in\N}\subseteq B$ such that $a=\bigvee_{n\in\N}b_n$. \end{definition}

Observe that in natural effect algebras of sets the family of the finite subsets of $\N$ is trivially a generator.

\begin{proposition}\label{BLv}
Let $L$ be an effect algebra, $G$ a normed group and $\mu:L\rightarrow G$ a $\sigma$-additive measure. Then $\mu$ is bounded on $L$ if and only if it is bounded on a generator of $L$.
\end{proposition}

\begin{proof}
Assume that $\|\mu(b)\|\leq M$ for all $b$ in a generator $B$ of $L$. Fix an arbitrary $a\in L$ and choose by hypothesis an increasing $(b_n)_{n\in\N}\subseteq B$ such that $a=\bigvee_{n\in\N}b_n$. In accordance with \cite[2.4]{Ab}, we have that $\left(\mu(b_n)\right)_{n\in\N}$ converges to $\mu\left(a\right)$. As a consequence, $\|\mu(a)\|\leq M$.
\end{proof}

We will define now a class of effect algebras which include the natural effect algebras of sets.

\begin{definition}
Let $L$ be an effect algebra. An orthogonal sequence $(b_n)_{n\in\N}$ of $L$ is called a natural basis of $L$ provided that $1=\bigoplus_{n\in\N}b_n$ and $b_n$ is minimal in $L\setminus\{0\}$ for every $n\in\N$.\end{definition}

Observe that in a natural effect algebra of sets a natural basis is the sequence $(\{n\})_{n\in\N}$.

\begin{proposition}\label{basisL}
Let $L$ be an effect algebra with the sRDP and $(b_n)_{n\in\N}$ a natural basis of $L$. For every $a\in L$ there exists a subsequence $(b_{n_k})_{k\in\N}$ of $(b_n)_{n\in\N}$ such that $a=\bigoplus_{k\in\N}b_{n_k}$. In particular, $B:=\left\{b_{i_1}\oplus\cdots\oplus b_{i_k}:i_1,\dots,i_k\in\N\right\}$ is a generator of $L$.\end{proposition}
\begin{proof}
Fix an arbitrary $a\in L$. Since $L$ is sRDP and $a\leq 1$, there exists an orthogonal sequence $(a_n)_{n\in\N}$ such that $a=\bigoplus_{n\in\N}a_n$ and $a_n\leq b_n$ for all $n\in\N$. By the minimality of the $b_n$'s we conclude that, for every $n\in\N$, either $a_n=b_n$ or $a_n=0$. This shows the existence of a subsequence $(b_{n_k})_{k\in\N}$ of $(b_n)_{n\in\N}$ such that $a=\bigoplus_{k\in\N}b_{n_k}$. Since $(b_{n_1}\oplus\cdots\oplus b_{n_k})_{k\in\N}$ is an increasing sequence of $B$, this also shows that $B$ is a generator of $L$.\end{proof}

\begin{lemma}\label{aaX}
Let $L$ be an effect algebra and $X$ a finite dimensional Banach space. If $\mu$ is a $\sigma$-additive $X$-valued measure and $(a_n)_{n\in\N}$ is an orthogonal sequence such that $\bigoplus_{n\in\N}$ exists, then $\sum_{n\in\N}\|\mu(a_n)\|<\infty$.
\end{lemma}
\begin{proof}
Observe that $\mu\left(\bigoplus_{n\in\N}a_n\right)=\sum_{n\in\N}\mu(a_n)$ is an unconditionally convergent series in a finite dimensional Banach space, therefore by the famous Riemann Series Theorem, the previous series is absolutely convergent.
\end{proof}

As a consequence of the previous lemma, $\sigma$-additive measures with values on a finite dimensional Banach space defined on $\sigma$-effect algebras are always absolutely additive.\\

Notice that the previous lemma cannot be generalized to infinite dimensional Banach spaces in virtue of the Dvoretzky-Rogers Theorem, which establishes the existence of non-absolutely convergent but unconditionally convergent series in every infinite dimensional Banach space.

\begin{theorem}
Let $L$ be an effect algebra with the sRDP and a natural basis $(b_n)_{n\in\N}$ and let $X$ be a finite dimensional Banach space. If $\mu$ is a $\sigma$-additive $X$-valued measure, then $\mu$ is bounded.\end{theorem}
\begin{proof}
In the first place notice that $M:=\sum_{n\in\N}\|\mu(b_n)\|<+\infty$ by Lemma \ref{aaX}. By Proposition \ref{basisL}, $B:= \left\{b_{i_1}\oplus\cdots\oplus b_{i_k}:i_1,\dots,i_k\in\N\right\}$ is a generator of $L$. According to Proposition \ref{BLv}(2), it suffices to show that $\mu$ is bounded on $B$. Let $b_{i_1}\oplus\cdots\oplus b_{i_k}\in B$. Then $$\left\|\mu\left(b_{i_1}\oplus\cdots\oplus b_{i_k}\right)\right\|=\left\|\sum_{j=1}^k\mu\left(b_{i_j}\right)\right\|\leq\sum_{j=1}^k\left\|\mu\left(b_{i_j}\right)\right\|\leq M.$$
Thus $\mu$ is bounded.\end{proof}

Examples of effect algebras with the sRDP and containing a natural basis are the Boolean algebras containing the finite subsets of $\N$.

The proof of the previous theorem can be adapted to version the previous theorem for measures valued in a normed group, but then we need those measures to be absolutely additive.

\begin{theorem}
Let $L$ be an effect algebra with the sRDP and a natural basis $(b_n)_{n\in\N}$ and let $G$ be a normed group. If $\mu$ is an absolutely additive $G$-valued measure, then $\mu$ is bounded.\end{theorem}
\begin{proof}
Since $(b_n)_{n\in\N}$ is orthogonal, by hypothesis we can consider $M:=\sum_{n\in\N}\|\mu(b_n)\|$. The rest follows as in the proof of the previous theorem.\end{proof}

The end of this section is on sequences of measures on effect algebras.

\begin{definition}
Let $L$ be an effect algebra, $G$ a topological (or normed if necessary) group and $(\mu_i)_{i\in\IN}$ a sequence of $G$-valued measures on $L$.
\begin{itemize}
\item $(\mu_i)_{i\in\IN}$ \emph{is pointwise bounded} if for every $a\in L$, there exists $\alpha_a>0$ such that $\|\mu_i(a)\|<\alpha_a$ for every $i\in\IN$.
\item $(\mu_i)_{i\in\IN}$ \emph{is uniformly bounded} if there exists $\alpha>0$ such that $\|\mu_i(a)\|<\alpha$ for every $a\in L$ and every $i\in\IN$.
\item $(\mu_i)_{i\in\IN}$ \emph{is uniformly exhaustive or uniformly strongly bounded} if for every orthogonal sequence $(a_n)_{n\in\IN}\subseteq L$, we have that $\lim_{n\in\IN}\mu_i(a_n)=0$ uniformly in $i\in\IN$.
\item $(\mu_i)_{i\in\IN}$ \emph{is uniformly strongly additive} if for every orthogonal sequence $(a_n)_{n\in\IN}\subseteq L$, we have that $\left(\sum_{k=1}^n\mu_i(a_k)\right)_{n\in\IN}$ converges uniformly in $i\in\IN$.
\item $(\mu_i)_{i\in\IN}$ \emph{is uniformly absolutely additive} if for every orthogonal sequence $(a_n)_{n\in\IN}\subseteq L$, we have that $\left(\sum_{k=1}^n\|\mu_i(a_k)\|\right)_{n\in\IN}$ converges uniformly in $i\in\IN$.
\item $(\mu_i)_{i\in\IN}$ \emph{is pointwise convergent} if there exists $\lim_{i\in\IN}\mu_i(a)$ for every $a\in L$.
\end{itemize}
\end{definition}

\begin{theorem}
Let $L$ be an effect algebra and ${(\mu_i)}_{i\in\IN}$ a sequence of $G$-valued measures on $L$ where $G$ is a complete Hausdorff topological group. Then:
\begin{enumerate}
\item ${(\mu_i)}_{i\in\IN}$ is uniformly strongly additive if and only if ${(\mu_i)}_{i\in\IN}$ is uniformly exhaustive.
\item If $G=\R$, then if ${(\mu_i)}_{i\in\IN}$ is uniformly exhaustive, then it is uniformly absolutely additive.
\end{enumerate}
\end{theorem}
\begin{proof}
\mbox{}
\begin{enumerate}
\item This follows from Theorem \ref{equiv}(1) applied to $\mu:= (\mu_i)_{i\in\N}:L\to G^{\N}$ where $G^{\N}$ is endowed with the topology of uniform convergence.
\item Suppose not.  Then there exist an orthogonal sequence ${(a_n)}_{n\in\IN}$ in $L$, $\vareps>0$, a strictly increasing $(n_r)_{r\in\N}$ in $\N$ and $i_r\in\N$ such that $\disp\sum_{n_{r-1}<k\leq n_{r}} |\mu_{i_r}(a_k)|\gq\vareps$.
Let  $\nu_n:=\mu_{i_n}$ for $n\in\IN$. Observe that ${(\nu_n)}_{n\in\IN}$ is uniformly exhaustive.
Let $A_r:=\{n_{r-1}+1,\dots, n_{r}\}$  and define
\begin{align*}
A_r^+& =\{k\in A_r:\nu_{n_{r-1}}(a_k)>0\},\\
A_r^-& =\{k\in A_r:\nu_{n_{r-1}}(a_k)<0\}.
\end{align*}
Consider now the sets
\begin{align*}
M_+& =\disp\biggl\{n\in\IN:\sum_{k\in A_n^+}\nu_n(a_k)\gq\frac{\vareps}{4}\biggr\},\\
M_-& =\disp\biggl\{n\in\IN:\sum_{k\in A_n^-} -\nu_n(a_k)\gq\frac{\vareps}{4}\biggr\}.
\end{align*}
Notice that $\IN=M_+\cup M_-$, and suppose, without lack of generality, that $M_+$ is infinite and put $M_+=\{m_1,m_2,m_3,\dots\}$ with $m_1<m_2<m_3<\dots$. For every $r\in\IN$, consider $$b_r=\disp\bigoplus_{k\in A_{m_r}^+}a_k\in L.$$ The $b_r$'s are in $L$ since ${(a_n)}_{n\in\IN}\subseteq L$ and ${(b_r)}_{r\in\IN}$ is an orthogonal sequence by definition. Observe that $\nu_{m_r}(b_r)=\disp\sum_{k\in A_{m_r}^+}\nu_{m_r}(a_k)\gq\frac{\vareps}{4}$ and denote ${(\xi_n)}_{n\in\IN}:={(\nu_{m_n})}_{n\in\IN}$, which is uniformly exhaustive. Thus, for $\disp\frac{\vareps}{4}>0$, there exists $n_0\in\IN$ such that if $n\gq n_0$ and $i\in\IN$, then $\disp|\xi_i(b_n)|<\frac{\vareps}{4}$. In particular, we have $|\xi_{n_0}(b_{n_0})|<\frac{\vareps}{4}$, but $$|\xi_{n_0}(b_{n_0})|=\nu_{m_{n_0}}(b_{n_0})=\disp\sum_{k\in A_{m_{n_0}}^+}\nu_{m_{n_0}}(a_k)\gq\frac{\vareps}{4},$$ which is a contradiction.
\end{enumerate}\end{proof}

\section{Vitali-Hahn-Saks, Nikodym and Grothendieck properties on effect algebras}

We will begin this section by introducing the Vitali-Hahn-Saks, Nikodym and Grothendieck properties on effect algebras.

\begin{definition}
We say that an effect algebra L has
\begin{itemize}
\item the Nikod\'ym property (for brevity N) if every sequence $(\mu_i)_{i\in\N}$ of bounded measures defined on $L$ with values in $\R$ which is pointwise bounded, is uniformly bounded;
\item the Vitaly-Hahn-Saks property (for brevity VHS) if every sequence $(\mu_i) _{i\in\N}$ of bounded measures defined on $L$ with values in $\R$ which is pointwise convergent, is uniformly exhaustive;
\item the Grothendieck property (for brevity G) if every sequence $(\mu_i) _{i\in\N}$ of bounded measures defined on $L$ with values in $\R$ which is pointwise convergent and uniformly bounded, is uniformly strongly additive.
\end{itemize}

If $L$ satisfies these conditions for sequences of $\sigma$-additive measures we will say that $L$ has the $\sigma$-Nikodym (for brevity $\sigma$-N), $\sigma$-Vitali-Hahn-Saks (for brevity $\sigma$-VHS) or $\sigma$-Grothendieck (for brevity $\sigma$-G) properties, respectively.\end{definition}

It is a classical result that a Boolean algebra has the VHS property if and only if it has the G and N properties (see \cite{Scha}). We will see that in the context of effect algebras, we achieve one of the implications, but if $L$ is an effect algebra with the VHS property, then it does not necessarily have the N property.

\begin{theorem}
Let $L$ be an effect algebra. If $L$ has the ($\sigma$-)G and ($\sigma$-)N properties, then $L$ has the ($\sigma$-)VHS property.
\end{theorem}
\begin{proof}
Let ${(\mu_i)}_{i\in\N}$ be a sequence of bounded real-valued ($\sigma$-additive) measures on $L$. If ${(\mu_i)}_{i\in\N}$ is pointwise convergent, then it is pointwise bounded and consequently it is uniformly bounded, since $L$ has the ($\sigma$-)N property. By hypothesis, $L$ has the ($\sigma$-)G property, so ${(\mu_i)}_{i\in\N}$ is uniformly strongly additive.
\end{proof}

We will show now an example of an effect algebra with the VHS property but failing to have the N property.

\begin{example}\label{VHSnoN}
Consider the effect algebra defined in Example \ref{ejemploEfectoNoAcotadaSumable} and define for each $i\in\IN$:
\begin{align*}
\mu_i(B_n)&=\left\{
\begin{array}{ll}
i , & \hbox{if $n=i$;}\\
0 , & \hbox{if $n\neq i$.}
\end{array}
\right.&
\mu_i(B_n^c)&=\left\{
\begin{array}{ll}
-i , & \hbox{if $n=i$;} \\
0, & \hbox{if $n\neq i$.}
\end{array}
\right. \\
\mu_i(\varnothing)&=0 & \mu_i(\IN)&=0.
\end{align*}
Observe that for each $i\in\IN$, $\mu_i$ is a bounded real-valued measure and the sequence $(\mu_i)_{i\in\N}$ is pointwise bounded, but it is not uniformly bounded, since $\mu_n(B_n)=n\to\infty$ as $n\to \infty$. This shows that $L$ does not have the N property. Observe that $L$ has the VHS property, since for every orthogonal sequence ${(a_j)}_{j\in\N}$, $\displaystyle\lim_{j\to\infty}\mu_i(a_j)=0$ uniformly in $i\in\IN$.
\end{example}

Observe that, by definition, if $L$ is an effect algebra with the VHS property, then it has the G property. Talagrand, in \cite{Tal84}, showed that a Boolean algebra with the G property does not necessarily have the N property by using the Continuum Hypothesis. Note that according to Example \ref{VHSnoN}, an effect algebra with the G property does not necessarily have the N property and we have not used the Continuum Hypothesis.\\

We will show next that in a certain class of effect algebras, which include the natural effect algebras of sets, the $\sigma$-VHS property, implies the $\sigma$-N property.

\begin{definition}
Let $L$ be an effect algebra. A subset $B$ of $L$ is said to be bounding provided that for every normed group $G$ and every pointwise bounded sequence $(\mu_k)_{k\in \N}$ of $\sigma$-additive $G$-valued measures on $L$ it is verified that $(\mu_k)_{k\in \N}$ is uniformly bounded on $L$ if and only if it is uniformly bounded on $B$.\end{definition}

The proof of Proposition \ref{BLv}(2) can be easily adapted to prove that generators are in fact bounding. In other words, we have the following lemma whose proof we omit.

\begin{lemma}\label{BLv2}
Let $L$ be an effect algebra, $G$ a normed group and $(\mu_k)_{k\in \N}$ a sequence of $\sigma$-additive $G$-valued measures on $L$. Then $(\mu_k)_{k\in \N}$ is uniformly bounded on $L$ if and only if it is uniformly bounded on a generator of $L$.\end{lemma}

Observe that in the previous lemma the sequence $(\mu_k)_{k\in \N}$ is not required to be pointwise bounded.

\begin{definition}
Let $L$ be an effect algebra. A natural subset $B$ of $L$ is a generator subset such that every orthogonal elements $b, c\in B$ verify that $b\oplus c\in B$ and for every $b\in B$ the set $B_b:=\{c\in B: c\perp b\}$ is bounding. We will say that $L$ is natural if it has a natural subset.\end{definition}

The following theorem generalizes \cite[Remark 3.2]{AiTa}.

\begin{theorem}
Let $L$ be a natural effect algebra. If $L$ has the $\sigma$-VHS property, then $L$ has the $\sigma$-N property.
\end{theorem}
\begin{proof}
Let $(\mu_i)_{i\in \N}$ be a sequence of real-valued bounded $\sigma$-additive measures on $L$ for which $(\mu_i(a))_{i\in\N}$ is bounded for every $a\in L$. We have to show that $(\mu_i)_{i\in \N}$ is uniformly bounded in $L$. According to Lemma \ref{BLv2}, it is sufficient to show that $(\mu_i)_{i\in \N}$ is uniformly bounded on $B$, where $B$ is any natural subset of $L$. Suppose to the contrary that $(\mu_i)_{i\in \N}$ is not uniformly bounded on $B$. We will follow an inductive process:
\begin{enumerate}
\item There exists $b_1\in B$ and $i_1\in\N$ such that $\left|\mu_{i_1}(b_1)\right|>1$. Let $B_1:=\{b\in B: b\perp b_1\}$. Notice that $(\mu_i)_{i>i_1}$ is not uniformly bounded in $B_1$ because $B_1$ is bounding.
\item There exists $b_2\in B_1$ and $i_2>i_1$ such that $\left|\mu_{i_2}(b_2)\right|>1$. Let $B_2:=\{b\in B: b\perp b_1\oplus b_2\}$. Notice that $(\mu_i)_{i>i_2}$ is not uniformly bounded in $B_2$ because $B_2$ is bounding.
\end{enumerate}
By induction we obtain an orthogonal sequence $(b_k)_{k\in\N}$ in $B$ and a subsequence $\left(\mu_{i_k}\right)_{k\in\N}$ of $(\mu_i)_{i\in \N}$ such that $\left|\mu_{i_k}(b_k)\right|>k$ for all $k\in \N$. Now $\left(\frac{1}{k}\mu_{i_k}\right)_{k\in\N}$ is a sequence of real-valued bounded $\sigma$-additive measures pointwise convergent to $0$ but $\left|\frac{1}{k}\mu_{i_k}(b_k)\right|>1$ for all $k\in \N$, which contradicts the fact that $L$ is $\sigma$-VHS.\end{proof}

The next result shows that the set of finite subsets of $\N$ is a natural subset of every natural effect algebra of sets.

\begin{theorem}
If $L$ is a natural effect algebra of sets, then $\phi_0(\N):=\{A\subset \N: A\text{ is finite}\}$ is a natural subset of $L$.
\end{theorem}
\begin{proof}
It is trivial that $\phi_0(\N)$ is a generator in $L$ and if $B,C\in \phi_0(\N)$ are disjoint, then $B\cup C\in \phi_0(\N)$. Fix an arbitrary $B\in\phi_0(\N)$ and let $\phi_B(\N):=\{C\in\phi_0(\N):B\cap C=\varnothing\}$. We have to show that $\phi_B(\N)$ is bounding. Indeed, let $G$ be a normed group and $(\mu_k)_{k\in\N}$ a pointwise bounded sequence of $G$-valued $\sigma$-additive measures on $L$ uniformly bounded on $\phi_B(\N)$. Let $M>0$ such that $\|\mu_k(C)\|\leq M$ for all $C\in \phi_B(\N)$ and all $k\in \N$. For every $n\in B$ let $M_n>0$ be so that $\left\|\mu_k(\{n\})\right\| \leq M_n$ for all $k\in \N$. We will finally prove that $(\mu_k)_{k\in\N}$ is uniformly bounded on $\phi_0(\N)$ by the bound $\sum_{n\in B}M_n + M$ (which implies that $(\mu_k)_{k\in\N}$ is uniformly bounded on $L$ because $\phi_0(\N)$ is a generator of $L$). Indeed, if $A\in\phi_0(\N)$, then
\begin{eqnarray*}
\|\mu_k(A)\|&\leq& \|\mu_k(A\cap B)\|+\|\mu_k(A\setminus B)\|\\
&\leq& \sum_{n\in A\cap B}\left\|\mu_k(\{n\})\right\|+ \|\mu_k(A\setminus B)\|\\
&\leq& \sum_{n\in B}\left\|\mu_k(\{n\})\right\|+ \|\mu_k(A\setminus B)\|\\
&\leq& \sum_{n\in B}M_n + M.
\end{eqnarray*}\end{proof}

We will conclude this section with a generalization of the famous Freniche's Lemma (see \cite{F}).

\begin{definition}\label{c_k}
An effect algebra $L$ is said to have
\begin{enumerate}
\item the orthogonally sequential property (OSP) if for every sequence $(b_k)_{k\in\N}$ of $L$ there exists an orthogonal sequence $(c_k)_{k\in\N}\subseteq L$ such that $c_k\leq b_k$ for all $k\in \N$ and if $a\in L$ and $k_0\in \N$ verify that $a\leq b_{k_0}$ and $b_i\leq 1\ominus a$ for all $i\in \N\setminus\{k_0\}$, then $a\leq c_{k_0}$;
\item the subsequential interpolation property (SIP) provided that for every orthogonal sequence $(a_n)_{n\in\N}\subseteq L$ and every infinite subset $M\subseteq \N$ there exists $a\in L$ and an infinite subset $N\subseteq M$ such that $a_n\leq a$ for all $n\in N$ and $a_n\leq 1\ominus a$ for all $n\in \N\setminus N$.
\end{enumerate}\end{definition}

Trivial examples of OSP effect algebras are Boolean algebras.

\begin{lemma}[Freniche's Lemma for effect algebras]
Let $G$ be a normed group. If $L$ is an OSP and SIP effect algebra, then for every $\varepsilon>0$, every orthogonal sequence $(a_n)_{n\in\N}\subseteq L$ and every exhaustive vector measure $\mu:L\rightarrow G$ there exist an infinite subset $M\subseteq \N$ and an element $a\in L$ such that $\|\mu(a)\|\leq \varepsilon$, $a_n\leq a$ for all $n\in M$ and $a\leq 1\ominus a_n$ for all $n\in \N\setminus M$.
\end{lemma}
\begin{proof}
Let $(M_k)_{k\in\N}$ be a pairwise disjoint sequence of infinite subsets of $\N$. By hypothesis, since $L$ has SIP, for every $k\in\N$ we can find $b_k\in L$ and an infinite subset $N_k\subseteq M_k$ such that $a_n\leq b_k$ if $n\in N_k$ and $b_k\leq 1\ominus a_n$ if $n\in\N\setminus N_k$. Let $(c_k)_{k\in\N}\subseteq L$ be as in Definition \ref{c_k}(1) for $(b_k)_{k\in\N}$. Since $\mu$ is exhaustive, $\disp\lim_{k\to\infty}\mu(c_k)=0$, therefore we can find $k_0\in \N$ such that $\left\|\mu\left(c_{k_0}\right)\right\|\leq \varepsilon$. Denote $a:=c_{k_0}$ and $M:=N_{k_0}$. Let us check that the theses of the lemma are verified:
\begin{itemize}
\item $\|\mu(a)\|=\left\|\mu\left(c_{k_0}\right)\right\|\leq \varepsilon$.
\item Fix an arbitrary $n\in M$. We know that $a_n\leq b_{k_0}$. On the other hand, if $i\in \N\setminus\{k_0\}$, since $N_{k_0}$ and $M_i$ are disjoint, we have that $n\in\N\setminus N_i$, thus $b_i\leq 1\ominus a_n$. Since $(c_k)_{k\in\N}$ verifies the properties of Definition \ref{c_k}, since $L$ has OSP we conclude that $a_n\leq c_{k_0}=a$.
\item Fix an arbitrary $n\in \N\setminus M$. We know that $b_{k_0}\leq 1\ominus a_n$ so $a:=c_{k_0}\leq b_{k_0}\leq 1\ominus a_n$.
\end{itemize}
\end{proof}

\section{The variation of a measure defined on an effect algebra}

We begin this final section by introducing the concept of variation of a measure on an effect algebra. We refer the reader to \cite{Bha,DiUhl77,DunSch88} for the classical concepts of variation of a measure on Boolean algebras of sets.

\begin{definition}
Let $L$ be an effect algebra, $G$ a normed group  and $\mu:L\rightarrow G$ a measure.
\begin{enumerate}
\item Let $e\in L$. We say that $\pi=\{e_1,e_2,\dots,e_n\}\subseteq L$ is a \emph{decomposition of $e$} if $\pi$ is an orthogonal set and $e=e_1\oplus\dots\oplus e_n$.
\item We define the \emph{variation of $\mu$ in $e\in L$} as the map $$\begin{array}{rrcl} |\mu|:&L&\to&[0,+\infty]\\ &e&\mapsto&|\mu|(e):=\disp \sup_{\pi\in\Pi}\sum_{e_i\in\pi}\|\mu(e_i)\|,\end{array}$$ where the supremum is taken over the family $\Pi$ of all decompositions $\pi=\{e_1,e_2,\dots,e_n\}$ of $e$.
\item We say that $\mu$ is \emph{of bounded variation} if $|\mu|(1)<+\infty$.
\end{enumerate}
\end{definition}

It is well-known that the variation of a measure defined on a Boolean algebra is also a measure (see \cite{DunSch88}). In the case of effect algebras, the sub-additivity of the variation is not granted (see Example \ref{ejemploNoSubad}).\\

The proof of the following lemma can  be easily seen by induction.

\begin{lemma}\label{rdpl}
Let $L$ be an RDP effect algebra. Consider $x\in L$ and $a_1,a_2,\dots,a_n\in L$ orthogonal such that $x\lq a_1\oplus a_2\oplus\dots a_n$. There exist $x_1,x_2,\dots,x_n\in L$ such that
\begin{align*}
x_1\lq a_1, x_2\lq a_2,\dots,x_n\lq a_n \\
x=x_1\oplus x_2\oplus\dots\oplus x_n.\\\end{align*}\end{lemma}

The next theorem shows that the variation of a measure defined on an effect algebra is always super-additive. The sub-additivity is obtained under the RDP hypothesis.

\begin{theorem}\label{VariacionMedida}
Let $L$ be an effect algebra, $G$ a normed  group and $\mu:L\rightarrow G$ a measure. Let $e,f\in L$ be orthogonal. Then:
\begin{enumerate}
\item $|\mu|(e)+|\mu|(f)\lq|\mu|(e\oplus f)$.
\item If $L$ is RDP, then $|\mu|(e\oplus f)\lq|\mu|(e)+|\mu|(f)$.
\end{enumerate}\end{theorem}
\begin{proof}
\mbox{}
\begin{enumerate}
\item The proof goes like in the Boolean case (see \cite{DiUhl77}).

\item If $L$ is RDP, given a decomposition $e\oplus f=a_1\oplus a_2\oplus\dots\oplus a_n$, since $e\lq e\oplus f= a_1\oplus a_2\oplus\dots\oplus a_n$, in virtue of Lemma \ref{rdpl}, there exist $e_1,e_2,\dots,e_n$ such that
\begin{align*}
e_1\lq a_1,e_2\lq a_2\dots e_n\lq a_n,\\
e=e_1\oplus e_2\oplus\dots\oplus e_n.
\end{align*}
Using the order defined on an effect algebra, there exist $f_1,\dots,f_n\in L$ such that $e_i\oplus f_i=a_i$ for all $i$'s. By the associative and commutative properties,
\begin{align*}
e\oplus f &=a_1\oplus a_2\oplus\dots\oplus a_n\\
&=(e_1\oplus f_1)\oplus(e_2\oplus f_2)\oplus\dots\oplus(e_n\oplus f_n)\\
&=(e_1\oplus e_2\oplus\dots\oplus e_n)\oplus(f_1\oplus f_2\oplus\dots\oplus f_n),
\end{align*}
and by the Cancellation Laws $f=f_1\oplus f_2\oplus\dots\oplus f_n$ and $f_i\lq f$ for all $i$'s. This shows that $\{e_i\}_1^n$ and $\{f_i\}_1^n$ are decompositions of $e$ and $f$, respectively. Thus,
\begin{align*}
\sum_{i=1}^n\|\mu(a_i)\|&=\sum_{i=1}^n\|\mu(e_i\oplus f_i)=\sum_{i=1}^n\|\mu(e_i)+\mu(f_i)\|\\
&\lq\sum_{i=1}^n\|\mu(e_i)\|+\sum_{i=1}^n\|\mu(f_i)\|\\
&\lq|\mu|(e)+|\mu|(f),
\end{align*}
Since this inequality is valid for every decomposition of $e\oplus f$, we conclude that
\begin{align}
|\mu|(e\oplus f)\lq|\mu|(e)+|\mu|(f).
\end{align}\end{enumerate}\end{proof}

\begin{corollary}\label{rdpm}
If $L$ is an RDP effect algebra, then the variation of a measure defined on $L$ is also a measure defined on $L$.
\end{corollary}

The following example shows that Corollary \ref{rdpm} does not hold true if we remove the hypothesis of RDP.

\begin{example}\label{ejemploNoSubad}
Consider $L=\{\varnothing,X^+,X^-,Y^+,Y^-,\IR^2\}$, where
\begin{align*}
X^+ &=\{(x,y)\in\IR^2:x\geq0\}\\
X^- &=\{(x,y)\in\IR^2:x<0\}\\
Y^+ &=\{(x,y)\in\IR^2:y\geq0\}\\
Y^- &=\{(x,y)\in\IR^2:y<0\}\\
\end{align*}
It is easy to show that $L$ is an effect algebra with $$A\oplus B:= A\cup B \Leftrightarrow A\cap B=\varnothing.$$ Define $\mu:L\rightarrow\IR$ such that
\begin{align*}
\mu(\varnothing)&= 0, & \mu(\IR^2)&=2, \\
\mu(X^+) &= 1, & \mu(X^-) &=1, \\
\mu(Y^+) &= 5, & \mu(Y^-) &=-3.\\
\end{align*}
It is easy to see that $\mu$ is a measure. However, observe that
\begin{align*}
|\mu|(X^+)&=1,\\
|\mu|(X^-)&=1,\\
|\mu|(X^+\cup X^-)&=|\mu|(\IR^2)=8.
\end{align*}
This proves that $|\mu|$ is not a measure.
\end{example}

The end of this section is devoted to show some properties of the variation of a measure defined on an effect algebra.

\begin{lemma}\label{Variacionlema}
Let $L$ be an effect algebra, $G$ a normed group and $\mu:L\rightarrow G$ a strongly additive measure. We have:
\begin{enumerate}
\item If $a\in L$, $\|\mu(a)\|\lq|\mu|(a)$.
\item If $a,b\in L$ are such that $a\lq b$, then $|\mu|(a)\lq|\mu|(b)$.
\end{enumerate}\end{lemma}
\begin{proof}
\mbox{}
\begin{enumerate}
\item Considering the trivial decomposition, the first item holds true.
\item Given a decomposition of $a$, we can add the element $b\ominus a$ and so we obtain a decomposition of $b$. By the definition of $|\mu|$ the thesis follows.
\end{enumerate}\end{proof}

\begin{proposition}\label{VariacionCotas}
Let $L$ be an effect algebra and $\mu:L\rightarrow\IC$ a bounded measure. For every $a\in L$,
\begin{enumerate}
\item $\disp\sup_{b\in L,\;b\lq a}|\mu(b)|\lq|\mu|(a).$
\item $\disp|\mu|(a)\lq 4 \sup_{b\in L,\;b\lq a}|\mu(b)|$.
\end{enumerate}
\end{proposition}
\begin{proof}
\mbox{}
\begin{enumerate}
\item The first item follows from the definition of $|\mu|$.
\item The proof goes like in the Boolean case (see \cite{DiUhl77}).
\end{enumerate}
\end{proof}

As a consequence of the previous proposition we obtain this simple but not less important corollary.

\begin{corollary}
Let $L$ be an effect algebra and $\mu:L\rightarrow\IC$ a  measure. Then, $\mu$ is a bounded measure if and only if $\mu$ is of bounded variation.
\end{corollary}

The next theorem generalizes for effect algebras an important result of Measure Theory on Boolean algebras (see \cite{DiUhl77}).

\begin{theorem}\label{VariacionFuertAcotada}
Let $L$ be an effect algebra and $\mu:L\rightarrow\IC$ a  bounded measure. Then $\mu$ is exhaustive if and only if $|\mu|$ is exhaustive.
\end{theorem}
\begin{proof}\mbox{}
\begin{enumerate}
\item[$\Leftarrow$] follows from  Lemma \ref{Variacionlema}(1).
\item[$\Rightarrow$] Suppose not. There exist $\vareps>0$ and an orthogonal sequence ${(a_n)}_{n\in \IN}\subseteq L$ such that for every $n\in\IN$, $|\mu|(a_n)\gq 4\vareps>0.$ By Proposition \ref{VariacionCotas}(2), for every $n\in\IN$ there exists $b_n\in L$ such that $b_n\lq a_n$ and $|\mu|(a_n)\lq 4|\mu(b_n)|$. This implies that if $n\in\IN$, $|\mu(b_n)|\gq\vareps$ and since ${(b_n)}_n$ is orthogonal, this contradicts the fact that $\mu$ is exhaustive.
\end{enumerate}
\end{proof}

\begin{proposition}\label{proposition}
Let $L$ be an effect algebra and $G$ a complete normed group. If $\mu:L\rightarrow G$ is a measure of bounded variation, then $\mu$ is strongly additive, hence $\mu$ is exhaustive by Theorem \ref{equiv}(1).
\end{proposition}
\begin{proof}
If ${(a_n)}_{n\in\IN}\subseteq L$ is an orthogonal sequence, then $$\sum_{n=1}^k\|\mu(a_n)\|\lq|\mu|\Biggl(\bigoplus_{n=1}^k a_n\Biggr){\lq}|\mu|(1).$$ Thus, $\disp\sum_{n\in\IN}\|\mu(a_n)\|\lq |\mu|(1)<\infty$, so $\disp\sum_{n\in\IN}\mu(a_n)$ is absolutely convergent and therefore convergent.
\end{proof}


\noindent \textbf{Acknowledgements.} The authors would like to express their deepest gratitude towards the referee for valuable comments and suggestions.

\end{document}